\documentclass[11pt]{article}
\usepackage{amssymb}
\usepackage{amsmath}
\usepackage{amsthm}
\usepackage[dvips]{graphicx}
\usepackage{enumitem}

\theoremstyle{plain}
\newtheorem{thm}{Theorem}

\newtheorem{prop}{Proposition}
\newtheorem{lem}{Lemma}

\newcommand{\R}{\mathbb{R}}

\newcommand{\N}{\mathbb N}
\newcommand{\Z}{\mathbb Z}

\newcommand{\vsp}{\vspace{0.1cm}}
\newcommand{\vs}{\vspace{0.3cm}}

\theoremstyle{remark}
\newtheorem{rem}{Remark}

\begin{document}

\title{On the sharp regularity for arbitrary actions of nilpotent groups on the interval:
the case of $N_4$}
\author{E. Jorquera, A. Navas \& C. Rivas}

\maketitle

\begin{abstract}
In this work, we determine the largest $\alpha$ for
which the nilpotent group of 4-by-4 triangular matrices with integer
coefficients and 1 in the diagonal embeds into the group of
$C^{1+\alpha}$ diffeomorphisms of the closed interval.

\end{abstract}

\vspace{0.5cm}

\section*{Introduction}

This work deals with the next general two-fold question:

\vsp\vsp

\textit{Given a group $G$ of orientation-preserving homeomorphisms of
a manifold $M$, is it conjugate to a group of diffeomorphisms of $M$?
If so, how smooth can this conjugate action be made?}

\vsp\vsp

In dimension larger than 1, the first half of the question has, in general, a negative
answer, even for the action of a single homeomorphism \cite{harrison}. However, in
the case where $M$ has dimension 1, this turns out to be very interesting, and
the answer deeply depends on the dynamical/algebraic structure of the
action/group considered. For instance, from the dynamical point of view, the
classical Denjoy theorem says that  a $C^2$ (more generally, $C^{1+bv}$)
orientation-preserving circle diffeomorphism with irrational rotation number
is necessarily conjugate to a rotation, hence minimal. On the other hand, in
lower regularity, there are the so-called Denjoy counterexamples, namely,
$C^{1+\alpha}$ diffeomorphisms with irrational rotation number that admit
wandering intervals; besides, every circle homeomorphism is conjugate to
a $C^1$ diffeomorphism. From the algebraic point of view, there is an
important obstruction for a group $G$ to admit a faithful action by
$C^1$-diffeomorphisms of a 1-manifold with boundary: every
finitely-generated subgroup of $G$ must admit a nontrivial
homomorphism onto $\Z$ (see \cite{thurston}; see also
\cite{navas-th} and \cite{BMNR}).

In this article, we focus on nilpotent group actions on the closed interval 
$[0,1]$. (Extensions of our results to the case of the circle are 
left to the reader.) The picture for Abelian group
actions is essentially completed by the works \cite{DKN,tsuboi}. For
non-Abelian nilpotent groups, an important theorem of J.Plante and
W.Thurston establishes that they do not embed in the group of $C^2$-diffeomorphisms
of $[0,1]$ (see \cite{plante-thurston}). However, according to B.Farb and J.Franks, every
finitely-generated, torsion-free nilpotent group can be realized as a group of
$C^1$ diffeomorphisms of $[0,1]$ (see also \cite{jorquera}). Motivated by 
this, we pursue the problem below, which was first addressed in \cite{FF} 
and stated this way in \cite{navas-book}. For the statement, 
recall that a diffeomorphism $f$ is said to be of class $C^{1+\alpha}$ if its
derivative is $\alpha$-Holder continous, that is, there exists $C>0$ such that \hspace{0.01cm}
$|f'(x)-f'(y)|\leq C|x-y|^\alpha$ \hspace{0.01cm} holds for all $x,y$. 

\vsp\vsp

\noindent{\bf Problem.} {\em
Given a nilpotent subgroup $G$ of $\mathrm{Homeo}_+([0,1])$, find the largest $\alpha$
such that $G$ embeds into the group $\mathrm{Diff}_+^{1+\alpha}([0,1])$
of $C^{1+\alpha}$ diffeomorphisms.}

\vsp\vsp

There are two results in this direction. First, in \cite{CJN} (see also \cite{navas-critic}),
the aforementioned Farb-Franks action of $N_d$, the nilpotent group of $d$-by-$d$
lower triangular matrices with integer entries and 1 in the diagonal, is studied in detail.
In particular, it is showed that this action cannot be made of class $C^{1+\alpha}$ for
$\alpha \geq \frac{2}{(d-1)(d-2)}$, yet it can be made $C^{1+\alpha}$ for any
smaller $\alpha$. Second, a recent result of K.Parkhe \cite{parkhe}
establishes that any action of a finitely-generated nilpotent group on $[0,1]$ is
topologically conjugate to an action by $C^{1+\alpha}$-diffeomorphisms for any
$\alpha < 1 / \kappa$, where $\kappa$ is the polynomial growth degree of the group.

For the particular case of $N_4$, the regularity obtained by Parkhe is hence smaller than
that of the Farb-Franks action, namely, $C^{1+\alpha}$ for $\alpha < 1/3$. Somehow
surprisingly, even this regularity is not sharp, as it is shown by our

\vsp\vsp

\noindent{\bf Theorem A.} {\em The
group $N_4$ embeds into $\mathrm{Diff}_+^{1+\alpha}([0,1])$ for every $\alpha<1/2$.}

\vsp\vsp

In \cite{CJN}, it is also shown that for  any $d\in \N$, there is a nilpotent group of nilpotence
degree $d$ embedded into $\mathrm{Diff}_+^{1+\alpha}([0,1])$, for any $\alpha < 1$. (This
is for instance the case of the Heissenberg group $N_3$.) This suggests that the optimal
regularity of a nilpotent group embedding into $\mathrm{Diff}_+([0,1])$ may not depend
on the degree of nilpotence. Our second result shows that, at least, this invariant is not
trivial, hence it is worth pursuing its study.

\vsp\vsp

\noindent{\bf Theorem B.} {\em The
group $N_4$ does not embed into $\mathrm{Diff}_+^{1+\alpha}([0,1])$ for any $\alpha>1/2$.}

\vsp\vsp

We point out that the $C^{3/2}$ regularity is not covered by our results, though we
strongly believe that $N_4$ does not admit an embedding in such regularity
(compare \cite{navas-critic}).

\vs

This article is organized as follows. In \S \ref{sec the group N_4}, we review some basic facts about the group $N_4$ such 
as normal forms. We also construct an action of $N_4$ on $\Z^3$ that preserves the lexicographic order on $\Z^3$; 
this action is inspired by the theory of left-orderable groups \cite{GOD}. In \S \ref{sec the embedding}, we show
that for any $\alpha<1/2$, the action of $N_4$ on $\Z^3$ can be projected into an action of $N_4$ on $[0,1]$ by $C^{1+\alpha}$
diffeomorphisms, which shows Theorem~A. Theorem~B in turn is proved in \S\ref{sec bounding the regularity}.

All actions considered in this work are by orientation-preserving maps.


\section{The group $N_4$}
\label{sec the group N_4}

Throughout this work, we use the following notation. Given two group
elements $x, y$, we let $[x,y] := xyx^{-1}y^{-1}$, and $x^y := yxy^{-1}$.
Recall that the derived series of a group $G$ is defined by $G^0 := G$
and $G^{i+1} := [G^i,G^i]$. The group $G$ is solvable of degree
$d$ if $G^d$ is trivial but $G^{d-1}$ is not. The central series of
$G$ is defined by $G^{(0)} := G$ and $G^{(i+1)} := [G,G^{(i)}]$.
The group $G$ is nilpotent of degree $\ell$ if $G^{(\ell)}$ is
trivial but $G^{(\ell-1)}$ is not.

The group $N_4$ is by definition the group of
matrices of the form
\begin{equation}\label{eq def N_4}\left(\begin{array}{cccc}1&0&0&0 \\e& 1&0&0\\a&f&1&0\\c&b&d&1
\end{array}\right),\end{equation}
where all the entries belong to $\Z$. We will use the
generating set $S$ of $N_4$ consisting of the matrices for which all
non-diagonal entries are $0$ except for one which is $1$. The elements of $S$
will be denoted by $e,f,d,a,b,c$, where each of these elements
represent the generating matrix with a 1 in the position
corresponding to the letter in (\ref{eq def N_4}); for example,
$$e=\left(\begin{array}{cccc}1&0&0&0 \\1& 1&0&0\\0&0&1&0\\0&0&0&1
\end{array}\right).$$

The reader can easily check that $N_4$ is isomorphic to the
(inner) semidirect product $\langle f,a,b,c\rangle \rtimes \langle
e,d\rangle$, where $\langle f,a,b,c\rangle\simeq \Z^4$ and $\langle
d,e \rangle\simeq \Z^2$. The conjugacy action of $\Z^2$ on $\Z^4$ is given by
\begin{equation}  \label{eq Ad_e} e: \hspace{0.1cm} f\mapsto fa^{-1}  \;,\; a\mapsto a \; ,\; b\mapsto bc^{-1} \; ,\; c\mapsto c,  \end{equation}
\begin{equation}  \label{eq Ad_d} d: \hspace{0.1cm} f\mapsto fb  \;,\; a\mapsto ac \; ,\; b\mapsto b \; ,\; c\mapsto c.  \end{equation}
In particular, $N_4$ is metabelian ({\em i.e.} it has solvability
degree 2). Further, $N_4$ has nilpotence degree 3: its lower
central series is given by
$$ N_4^{(1)}=\langle a,b,c\rangle  \, ,\; N_4^{(2)}=\langle c \rangle \, , \; N_4^{(3)} =\{id\} .$$

It follows from equations (\ref{eq Ad_e}) and (\ref{eq Ad_d}) that  any element of $N_4$
can be written in a unique way as $$f^{n_1} e^{n_2} d^{n_3} a^{n_4} b^{n_5} c^{n_6},$$
where the exponents $n_i$ belong to $\Z$. This will be our preferred normal form. It allows proving the next

\vsp

\begin{lem} \label{lem inyectivo} Let $\phi: N_4\to G$ be a group homomorphism such that
$\phi(c)$ is a nontrivial element of $G$ with infinite order. Then $\phi$ is an embedding.
\end{lem}

\noindent{\em Proof:}  We first observe that, for $(n_1,n_2) \neq (0,0)$,
\begin{equation} \label{eq referee}[\phi(d^{n_1}e^{n_2}), \phi(a^{n_1}b^{-n_2}c^{n_3})]
= \phi([d^{n_1}e^{n_2}, a^{n_1}b^{-n_2}c^{n_3}])= \phi(c^{n_1^2+n_2^2}).\end{equation}
By the hypothesis, $\phi(c^{n_1^2 + n_2^2}) \neq id$, which
implies that the restriction of $\phi$ to both $\langle a,b,c\rangle$
and $\langle d,e \rangle$ is an embedding.

Further, for $(n_1,n_2)\neq (0,0)$, we have
$$\phi([d^{n_1}e^{n_2} a^{n_3}b^{n_4}c^{n_5},a^{n_1}b^{-n_2}])=
\phi([d^{n_1}e^{n_2},a^{n_1}b^{-n_2}])=\phi( c^{n_1^2+n_2^2})\neq id, $$
thus the restriction of $\phi$ to $\langle d,e,a,b,c\rangle$ is an embedding. Finally we have that, for  $n_0\neq 0$,
$$\phi([f^{n_0}e^{n_1}d^{n_2} a^{n_3}b^{n_4}c^{n_5}, e])=\phi(a^{n_0}c^{n_4})\neq id.$$
Hence, $\phi$ is injective.  $\hfill\square$

\vs

\begin{rem}\label{rem one orbit} An immediate consequence of Lemma \ref{lem inyectivo} is that for every
faithful action of $N_4$ by homeomorphisms of $[0,1]$, there is a point $x_0\in (0,1)$ such that $N_4$
acts faithfully on its orbit. Indeed, it suffices to consider $x_0$ as being any point moved by $c$.
\end{rem}

\vsp

We next construct an action of $N_4$ by homeomorphisms of  $[0,1]$.
Our method is close to the construction of Farb and Franks, who first built an action of
$N_4$ on $\Z^3$ and then project it to an action on $[0,1]$; see \cite{FF} or \cite{CJN}.
However, it should be emphasized that our action is different, which allows improving 
the regularity. We begin with

\vsp

\begin{prop} \label{prop action on Z^3} Let $\tilde e$, $\tilde f$, $\tilde d$ ,$\tilde a$, $\tilde b$, and $\tilde c$ 
be the maps from $\Z^3$ to $\Z^3$ defined by:
\begin{eqnarray}
\notag \tilde e :(i,j,k)\mapsto (i+1,j,k), \\
\notag \tilde d :(i,j,k)\mapsto (i,j+1,k),\\
\notag \tilde f :(i,j,k)\mapsto(i,j,k-ij),\\
a':(i,j,k)\mapsto (i,j,k-j),\\
\notag \tilde b :(i,j,k)\mapsto(i,j,k+i),\\
\notag \tilde c :(i,j,k)\mapsto (i,j,k+1).
 \end{eqnarray}
Then the group $\tilde N$ generated by $\langle \tilde e ,\tilde f ,\tilde d ,\tilde a ,\tilde b ,\tilde c \rangle$ is isomorphic to $N_4$.
\end{prop}

\noindent{\em Proof:} It follows from the definition that
$\tilde f,\tilde a,\tilde b$ and $\tilde c$ commute, and that the subgroup of $\tilde N$ 
that they generate is normal and isomorphic to $\Z^4$. Further, the subgroup
generated by $\{\tilde e,\tilde d\} $ is Abelian, and its action by conjugation
on $\langle \tilde f,\tilde a,\tilde b,\tilde c \rangle$ mimics equations (\ref{eq Ad_e})
and (\ref{eq Ad_d}). Therefore, by Lemma \ref{lem inyectivo},
the application $x\mapsto \tilde x$, with $x\in\{e,d,f,a,b,c\}$, induces
an isomorphism between $N_4$ and $\tilde N$. $\hfill\square$

\vsp

We now let $\left( I_{i,j,k} \right)_{(i,j,k)\in\Z^3}$ be a family of disjoint open intervals disposed on $[0,1]$ respecting
the (direct) lexicographic order of $\Z^3$, that is, $I_{i,j,k}$ is to the left of $I_{i',j',k'}$ if and only if $(i,j,k) \prec (i',j',k')$,
where $\preceq$ is the lexicographic order on $\Z^3$. Assume further that the union of this family of intervals is dense in
$[0,1]$.  Then, by some abuse of notation, we can define $e,d,f$ to be the unique homeomorphism of
$[0,1]$ whose restriction to each of the intervals $I_{i,j,k}$ is affine and send, respectively,
\begin{eqnarray}
\notag e: I_{i,j,k}\mapsto I_{i+1,j,k}, \\
\label{eq accion afin} 
d: I_{i,j,k }\mapsto I_{i,j+1,k},\\
\notag f: I_{i,j,k}\mapsto I_{i,j,k+ij}.
\end{eqnarray}
Since an affine map fixing a bounded interval must be the identity, Proposition \ref{prop action on Z^3} implies that the 
homeomorphisms $e,d,f$ generate a subgroup of $\mathrm{Homeo}_+([0,1])$ isomorphic to $N_4$. In  order to show 
Theorem A, in \S \ref{sec the embedding}, we will use, instead of affine maps, the so-called Pixton-Tsuboi family of 
local diffeomorphisms \cite{pixton,tsuboi}.

\vsp

\begin{rem} At first glance, this action may look strange. However, it naturally appears when considering total order 
relations that are invariant under left-multiplication ({\em left-orders}, for short); see \cite{GOD}. Namely, we may first 
endow the subgroups $\langle e,d \rangle \sim \mathbb{Z}^2$ and $\langle f,a,b,c \rangle \sim \mathbb{Z}^4$ 
with the left-orders $\preceq_1$ and $\preceq_2$, respectively, and then consider the convex extension of these,  
which is a left-order $\preceq$ on $N_4$.  In our construction, on the one hand, we let $\preceq_1$ 
be the lexicographic order for which $e$ is (infinitely) larger than $d$. On the other hand, we let $\preceq_2$ be the lexicographic order in which $c$ is the largest generator. Proceeding this way, the {\em dynamical realization} of the order $\preceq$ is 
an action of $N_4$ on the real line that is semiconjugated into the action above.   
\end{rem}

\vsp


\section{Bounding the regularity}
\label{sec bounding the regularity}


In this section, we show that the group $N_4$ does not embed in
$\mathrm{Diff}_+^{1+\alpha}([0,1])$ provided that $\alpha>1/2$. We first
reduce Theorem B to a combinatorial statement, namely Lemma \ref{lem main} below.


\subsection{The combinatorics prevents an embedding}

To state the main combinatorial lemma (whose proof is postponed to \S \ref{sec proof of main lemma}.), we introduce a notation that will be used throughout all \S2. 

\vsp

Given a nilpotent group $G$ acting by homeomorphisms of $[0,1]$, a point $x_0\in [0,1]$, and an element $g\in G$, we define
\begin{equation}\notag \label{eq intervalos}J_g(x_0):=[\inf_n g^n(x_0), \sup_n g^n(x_0)].\end{equation}
Since $G$ is nilpotent, given any $h\in G$, we have that the intervals $h(J_g)$ and $J_g$ either are equal or have disjoint interior (otherwise, one can build a free subsemigroup inside $G$; see for instance \cite[\S3.2]{GOD}). In the latter case, we will say that $h$ {\em moves} $J_g$.  We have

\begin{lem}
\label{lem main} Suppose that $N_4$ is faithfully acting on $[0,1]$ by $C^{1+\alpha}$-diffeomorphisms for some $\alpha>1/2$. 
Then there exist $g_1$, $g_2$, $g_3$ in $N_4$ and $x_0\in[0,1]$ such that: 
\begin{enumerate}
\item  $J_{g_3}(x_0)$ is not reduced to a point.
\item The element $g_2$ moves $J_{g_3}(x_0)$ and the element $g_1$ moves $J_{g_2}(x_0)$.
\item The elements $g_1,$ $g_2$, and $g_3$  pairwise commute. In particular, the subgroup $\langle g_1,g_2,g_3\rangle$ is isomorphic to $ \Z^3$.
\end{enumerate}
\end{lem}

\vsp

Lemma \ref{lem main} provides us enough combinatorial information about the dynamics of $N_4$ to prove Theorem B. In concrete terms, looking for a contradiction suppose that $N_4$ is faithfully acting by $C^{1+\alpha}$-diffeomorphisms for some $\alpha>1/2$, and let  $g_1,\; g_2,\; g_3$ and $x_0$ be the {\em elements} provided by the conclusion of Lemma \ref{lem main}. Then the only element in the Abelian group $\langle g_1, g_2, g_3\rangle$ fixing $x_0$ is the trivial one.
Further, by eventually changing some of $g_1, g_2,g_3$ by their inverses, we can suppose that they all move $x_0$ to the right. Hence,
if we define $I_{0,0,0}$ as the interval $(x_0,g_3(x_0))$ and $I_{n_1,n_2,n_3} := g_1^{n_1} \, g_2^{n_2} \, g_3^{n_3} \,(I_{0,0,0})$,
then the intervals $I_{i,j,k}$ are pairwise disjoint, they are disposed on $[0,1]$ respecting the lexicographic order of the indices, and
$$g_1 (I_{i,j,k})= I_{i+1,j,k}\, ,\; g_2(I_{i,j,k})=I_{i,j+1,k}\, , \; g_3(I_{i,j,k}) = I_{i,j,k+1}.$$
A contradiction is then provided by the following theorem from \cite{navas-critic} (see also \cite{DKN})

\vsp

\begin{thm} \label{thm navas}\textit{Let $k\!\geq\!2$ be an integer, and let $f_1,\ldots,f_{k}$
be commuting $C^1$-diffeomorphisms of $[0,1]$.
Suppose that there exist disjoint open intervals
$I_{n_1,\ldots,n_{k}}$ disposed on $(0,1)$ respecting the lexicographic order and so
that for all $(n_1,\ldots,n_{k}) \! \in \! \mathbb{Z}^{k}$ and all $i\! \in\! \{1,\ldots,k\}$,
$$f_i(I_{n_1,\ldots,n_i,\ldots,n_{k}}) = I_{n_1,\ldots,n_i+1,\ldots,n_{k}}.$$
Then $f_1,\ldots,f_{k-1}$ cannot be all simultaneously of class
$C^{1+1\!/\!(k-1)}$ provided that $f_{k}$ is of class
$C^{1+\alpha}$ for some $\alpha>0$.}
\end{thm}


\subsection{Proof of Lemma \ref{lem main}}
\label{sec proof of main lemma}

As discussed in the previous section, in order to finish the proof of Theorem B, we need to prove Lemma \ref{lem main}.
A first crucial step is given by the next result, which can be thought of as a version of Denjoy's theorem on the interval
and corresponds to an extension of \cite[Theorem C]{DKN} for the case where the maps are not assumed to commute.

\vsp

\begin{thm}
\label{thm non minimal} Given an integer $d\geq 2 $ and $\alpha>1/d$, suppose that $G$ is a subgroup of
$\mathrm{Diff}_+^{1+\alpha}([0,1])$ whose action is semiconjugated to a free action by translations of $\Z^d$.
Then $G$ acts minimally on $(0,1)$, and it is hence Abelian.
\end{thm}




\vsp

\noindent{\em Proof:}  Looking for a contradiction, we suppose that
the action of $G$ is not minimal. We let $I$ be a maximal open interval that is mapped into a single 
point by the semiconjugacy into a group of translations, and we let $f_1,\ldots,f_d \in G$ be elements 
whose semiconjugate images generate $\mathbb{Z}^d$. Changing the $f_i$'s by their inverses if 
necessary, we may assume that they all move points inside $(0,1)$ to the left.

We follow the proof of \cite[Theorem 4.1.37]{navas-book}, where the $f_i$'s are assume to commute. 
Although in our situation the $f_i$'s do not {\em a priori} commute, they do commute on the closure 
$\Lambda$ of the orbit of the endpoints of $I$.  This allows applying all arguments of \cite{navas-book} except the 
last one, provided we consider the underlying Markov process directly on intervals. More precisely, assume that 
all the $f_i$'s are tangent to the identity at the origin (the other case works almost verbatim to \cite{navas-book}; 
alternatively, use the M\"uller-Tsuboi trick \cite{tsuboi} to ensure flatness). 
Then consider the Markov process on $\mathbb{N}_0^d$ with transition probabilities
$$p \big((n_1,\ldots,n_i,\ldots,n_d) \to (n_1,\ldots,1+n_i,\ldots,n_d) \big)
:= \frac{1+n_i}{d + n_1 + \ldots + n_d}.$$
Denote by $\Omega$ the space of infinite paths $\omega$
endowed with the induced probability measure $\mathbb{P}$. Let
$S \!: \Omega \rightarrow \mathbb{R}$ be defined by
$$S (w) = \sum_{k \geq 0} |I_{\omega_k}|^{\alpha},$$
where $w_k = (n_{1,k},\ldots,n_{d,k}) \in \mathbb{N}_0^d$ denotes the position of $w$ at
time $k$, and $I_{n_1,\ldots,n_d} := f_1^{n_1}\ldots f_d^{n_d} (I)$.
Since $\alpha > 1/d$, this function has a finite expectation (see \cite{DKN}). Thus, its value at
a generic random sequence $\omega$ is finite. As in the proof of  \cite[Theorem 4.1.37]{navas-book}, 
if for such a sequence we denote $h_k := f_1^{n_{1,k}} \cdots f_{d}^{n_{d,k}}$, then we have 
\begin{equation}\label{just-this}
\frac{D  h_k (y)}{D  h_k (x)} \leq C
\end{equation}
for all $k \geq 1$ and all $x,y$ in $\bar{I} \cup \bar{J}$, where $C$ only depends on $\omega$ and the $\alpha$-H\"older 
constants of the derivatives of the 
$f_i$'s, and $J$ is any interval that is next to $I$ and has length smaller than $| I | / C$. By 
the maximality of $I$, there must exist some $h \in G$ mapping  $I$ into $J$. We then notice that, 
if $J$ has endpoints in $\Lambda$, then for all $k \geq 1$ we have 
$$\frac{| h(I) |}{| I |} 
= \frac{| h_k^{-1} h h_k (I) |}{| I |} 
= \frac{|h_k(I)|}{| I |} \cdot \frac{|hh_k(I)|}{| h_k(I) |} \cdot \frac{| h_k^{-1} h h_k (I) |}{| h h_k (I)|}.$$
In the product above, the middle quotient converges to $Dh (0) = 1$ as $k$ goes to infinite. Besides, 
the first and the third quotients are respectively equal to $Dh_k (x_k)$ and $1/ D h_k (y_k)$ for certain  
points $x_k \in \bar{I}$ and $y_k \in \bar{J}$. Using (\ref{just-this}), we conclude that $| h(I )| / | I | \geq 1/C$. 
However, this is impossible if $J$ was chosen small-enough so that $|J| < 1/C$. $\hfill\square$

\vsp\vsp\vsp

To finish the proof of Lemma \ref{lem main}, recall that every finitely-generated nilpotent group $G$ of homeomorphisms of $(0,1)$  preserves a
nontrivial Radon measure $\mu$ on $(0,1)$; see \cite{plante} or \cite{navas-book}. This measure induces a group homomorphism,
the so-called {\em translation number homomorphism} $\tau_\mu \!: G\to \R$, whose kernel coincides with the set of elements in $G$ 
having fixed points, and every such element must fix all points in $supp(\mu)$, the support of $\mu$. Moreover, if $\tau_\mu(G)$ has rank 
greater than or equal to 2, then $G$ is semiconjugate to the group of translations $\tau_\mu(G)$. In particular, from this we obtain

\vsp\vsp

\begin{lem} \label{lem heisemberg}Suppose the Heisemberg group $N_3\simeq \langle h_1. h_2, h_3\mid [h_1,h_2]=h_3, \;\;h_ih_3=h_3 h_i \;(i=1,2)\rangle$ 
is faithfully acting by homeomorphisms of $[0,1]$. If $x$ is not fixed by $h_3$, then at least one of $h_1, h_2$ moves $J_{h_3}(x)$.
\end{lem}

\vsp\vsp



We are now in position to give the 

\vsp\vsp

\noindent{\em Proof of Lemma \ref{lem main}:} Suppose $N_4$ faithfully acts by
$C^{1+\alpha}$ diffeomorphisms of $[0,1]$ for some $\alpha >1/2$. We let $x_0$ be a point moved by $c$. By Remark \ref{rem one
orbit},  $N_4$ faithfully acts on its orbit. To simplify the notation, for $g\in N_4$, the interval $J_g(x_0)$ will be denoted by $J_g$. 

\vsp

The key observation is that in $N_4$ there are many isomorphic copies of the Heissemberg group $N_3$ so there are many instances in which we can apply Lemma \ref{lem heisemberg}. The reader can easily check that, for example, the subgroups $$\langle e,b,c\rangle, \; \langle d,a,c\rangle,\; \langle f,d,b\rangle, \; \langle f,e,a\rangle$$ are all isomorphic to $N_3$ (the right-most generator being the generator of the center of $N_3$).

We let $g_3 := c$. Since $J_c$ is not reduced to a point, the first part of the conclusion of Lemma \ref{lem main} is satisfied. 
In order to find $g_2$ and $g_3$, we distinguish two cases:

\vspace{0.15cm}

{\bf Case 1:} Either $a$ or $b$ moves $J_c$.

\vsp

Suppose $a$ moves $J_c$. Then from Lemma \ref{lem heisemberg} applied to $\langle f,e,a \rangle$
we have that either $f$ or $e$ moves $J_a$. Then we can let $g_2:=a$ and $g_1$ be an element in $\{f,e\}$ that moves $J_a$. For these elements the conclusion of Lemma \ref{lem main} holds.

The case where $b$ moves $J_c$ works in the same way but looking at $\langle f,d,b \rangle$ instead of $\langle f,e,a \rangle$.

\vspace{0.15cm}

{\bf Case 2:} Both $a$ and $b$ fix $J_c$.

\vsp

Consider the group $\langle e,d\rangle\simeq \Z^2$ acting on the smallest possible interval containing $x_0$, that is, the convex closure $I$ 
of the $\langle e,d \rangle$-orbit of $x_0$. Observe that  $I$ is not contained in $J_c$ since in that case both $d$ and $a$  would fix $J_c$, thus 
contradicting Lemma \ref{lem heisemberg} applied to $\langle d,a,c \rangle$. In particular, the action of $\langle e,d \rangle$ on $I$ is not minimal. Theorem \ref{thm non minimal} then implies that the action of $\langle e,d\rangle $ on $I$ is not semi-conjugated to an action by translation of $\Z^2$, so there must be $h_0\in \langle e,d\rangle$ with translation number (over $I$) equal to zero.

If $h_0$ moves $J_c$ we are done, since we can let $g_2 :=h_0$ and $g_3$ be any element in $\langle e,d \rangle $ with non-trivial translation number. We claim that this is always the case; more precisely, we claim that any $h\in\langle e,d \rangle$ different from the identity moves $J_c$. Indeed, if $h=e^nd^m$ fixes $J_c$, then the group $H=\langle e^n d^m, a^m b^{-n}\rangle$ fixes $J_c$. But, if $(n,m)\neq (0,0)$, then equation (\ref{eq referee}) implies that $H$ is isomorphic to the Heisemberg group $N_3$ with center generated by $c^{n^2+m^2}$. A contradiction is then provided by Lemma \ref{lem heisemberg} and the fact that $J_c=J_{c^k}$ for any $k\neq 0$. 

\vsp

This finishes the proof of the Lemma \ref{lem main}, and hence that of Theorem~B.~$\hfill\square$


\section{The embedding}
\label{sec the embedding}


We next prove Theorem A. For the rest of this work, we fix $\alpha$ such that $0<\alpha<1/2$. In order to
produce an embedding of $N_4$ into $\mathrm{Diff}_+^{1+\alpha}([0,1])$, we will project to the interval the 
action provided by Proposition \ref{prop action on Z^3} using the so-called Pixton-Tsuboi maps \cite{pixton,tsuboi}. 
This technique is summarized in the next 

\vsp

\begin{lem} \label{pixton} 
There exists a family of $C^{\infty}$ diffeomorphisms $\varphi_{I',I}^{J',J}: I \to J$ between intervals $I,J$,
where $I'$ (resp. $J'$) is an interval contiguous to $I$ (resp. $J$) by the left, such that:
\begin{enumerate}
\item
(Equivariance) For all $I,I',J,J',K,K'$ as above,
$$\varphi_{J',J}^{K',K} \circ \varphi_{I',I}^{J',J} = \varphi_{I',I}^{K',K}.$$

\item
(Derivatives at the endpoints) For all $I,I',J,J'$, 
$$D \varphi_{I,I'}^{J,J'} (x_{-}) = \frac{|J'|}{|I'|}, \quad D \varphi_{I,I'}^{J,J'} (x_{+}) = \frac{|J|}{|I|},$$
where $x_{-}$ (resp. $x_+$) is the left (resp. right) endpoint of $I$. 

\item
(Regularity) There is a constant $M$ such that for all $x \in I$, we have 
$$D \log (D \varphi_{I',I}^{J',J}) (x) \leq \frac{M}{|I|} \cdot \left| \frac{|I|}{|J|} \frac{|J'|}{|I'|} - 1\right|$$
provided that \hspace{0.01cm} $\max \{ |I'||I|,|J'|,|J| \} \leq 2 \min \{ |I'||I|,|J'|,|J| \}$. 
\end{enumerate}
\end{lem}

\vsp\vsp

To produce our action, we let $I_{i,j,k}$ be a collection of intervals indexed by $\Z^3$ whose union is dense in $[0,1]$ 
and that are disposed preserving the lexicographic order of $\Z^3$. We then define the homeomorphisms $d,e,f$ of 
$[0,1]$ as those whose restrictions to $I_{i,j,k}$ coincide, respectively, with  
$$\varphi_{I_{i,j,k-1},I_{i,j,k}}^{I_{i+1,j,k-1},I_{i+1,j,k}}, \quad 
\varphi_{I_{i,j,k-1},I_{i,j,k}}^{I_{i,j+1,k-1},I_{i,j+1,k}}, \quad \mbox{and} \quad 
\varphi_{I_{i,j,k-1},I_{i,j,k}}^{I_{i,j,k+ij-1},I_{i,j,k+ij}}.$$
By ({\em Equivariance}), this produces a faithful action of $N_4$ by homeomorphisms of $[0,1]$.

\vsp

\begin{prop} For an appropriate choice of the lengths $|I_{i,j,k}|$, the homeomorphisms $e,f,d$ 
are simultaneously of class $C^{1+\alpha}$.
\end{prop}

\vsp

The rest of this work is devoted to the proof of this result. To begin with, 
we let $p,q,r$ be positive reals for which the following conditions hold:\\
(i)   $\alpha + r \leq 2$,\\
(ii)  $4r\leq p$,\\
(iii) $4r\leq q$,\\
(iv)   $4\leq p(1-\alpha)$,\\
(v)  $4\leq q(1-\alpha)$,\\
(vi) $1/p+1/q+1/r<1$,\\
(vii) $\alpha\leq \frac{1}{p}+\frac{1}{r}$ \hspace{0.1cm} and \hspace{0.1cm} $\alpha\leq \frac{r}{p(r-1)}$,\\
(viii)  $\alpha\leq \frac{1}{q}+\frac{1}{r}$ \hspace{0.1cm} and \hspace{0.1cm} $\alpha\leq \frac{r}{q(r-1)}$.\\
For example, we can take $p=q:=4/\alpha$ and $r:=4/3$.

Now, let $I_{i,j,k}$ be an interval such that
$$\vert I_{i,j,k}\vert:=\frac{1}{\vert i \vert^{p}+\vert j\vert^{q}+\vert k\vert^{r}+1}.$$
Condition (vi) ensures that 
$$\sum_{(i,j,k) \in \Z^3} \big| I_{i,j,k} \big| < \infty,$$
hence the $I_{i,j,k}$'s can be disposed on a finite interval respecting the lexicographic order. 
This interval can be though of as $[0,1]$ after renormalization. 

Observe also that conditions (i) to (viii) can only be satisfied for $\alpha < 1/2$. Indeed, the second part  
of conditions (vii) and (viii) together imply that 
$$2\alpha \big( 1-\frac{1}{r} \big) \leq \frac{1}{p}+\frac{1}{q}.$$ 
Then using (vi), one easily concludes that $2\alpha<1$, that is, $\alpha < 1/2$.

\vsp

It is proved in \cite{CJN} that, with any choice of lengths as above, the maps $e$ and $d$ are $C^{1+\alpha}$ 
diffeomorphisms. More precisely, in \cite[\S3.3]{CJN} it is shown that, under condition (vii), the diffeomorphism 
$e$ is of class $C^{1+\alpha}$. Indeed, the second half of condition (vii) corresponds to condition 
($\mathrm{iii}_{\mathrm{B}}$) in \cite{CJN}, while the first half of condition (vii), although not explicitly stated in 
\cite{CJN}, corresponds to the right form of condition ($\mathrm{v}_{\mathrm{B}}$) therein to show the regularity 
of $e$; see, for instance, \cite[page 125, line 9]{CJN}. 

An analog argument applies to $d$ under condition (viii). In order to conclude, below we will develop several 
slight modifications of some of the arguments of \cite{CJN} to show the next lemma, which closes the proof 
of Theorem A.\vsp

\begin{lem} For any choice of lengths of intervals satisfying properties {\em (i),...,(viii)} above, 
the homeomorphism $f$ is a $C^{1+\alpha}$ diffeomorphism.
\end{lem}

\vsp

Notice that this lemma is equivalent to that the expression 
$$\frac{\big| \log Df(x) - \log Df(y) \big|}{|x-y|^{\alpha}}$$
is uniformly bounded (independently of $x$ and $y$).  
To check this, due to property ({\em Derivatives at the endpoints}) above, it suffices 
to consider points $x,y$ in intervals $I_{i,j,k}$ and $I_{i,j,k'}$, respectively; this means 
that the first ``two levels'' $i$ and $j$ coincide (compare \cite[\S 3.3, III]{CJN}).  We will first deal with 
the case where the points $x,y$ belong to the same interval $I_{i,j,k}$, and then with that where these  
points lie in intervals of this form but with different indices $k,k'$. 

\vspace{0.35cm}

{\bf Case 1:} The points $x,y$ belong to the same interval $I := I_{i,j,k}$.

\vsp

In this case, as $|x-y| \leq |I|$, from ({\em Regularity}) in Lemma \ref{pixton} and the Mean Value Theorem  
we deduce that we need to find an upper bound for 
$$ \frac{1}{|I|^\alpha} \left| \frac{|I| |J'|}{|I'||J|} -1 \right|,$$
where $J$ denotes $f(I) := I_{i,j,k+ij}$, $I' := I_{i,j,k-1}$, and $J' := I_{i,j,k+ij-1}$.

\vs

{\bf Case 2:} The points $x,y$ lie in different intervals, say $x \in I_{i,j,k}$ and $y \in I_{i,j,k'}$, with $k' > k$.

\vsp

Here, using \cite{tsuboi}
(more precisely, 
\cite[(20)]{CJN}), it readily follows from the triangular inequality that $\vert \log Df (x) - \log Df (y) \vert$ is smaller 
than or equal to 
\begin{small}
$$\left \vert \log \!\frac{\vert I_{i,j,k+ij}\vert}{\vert I_{i,j,k}\vert} \!-\!
\log \!\frac{\vert I_{i,j,k'+ij} \vert}{\vert I_{i,j,k'}\vert}\right \vert +
\left \vert \log \!\frac{\vert I_{i,j,k+ij-1}\vert}{\vert I_{i,j,k-1}\vert} \!-\!
\log \!\frac{\vert I_{i,j,k+ij} \vert}{\vert I_{i,j,k}\vert}\right \vert +
\left \vert \log \!\frac{\vert I_{i,j,k'+ij-1}\vert}{\vert I_{i,j,k'-1}\vert} \!-\!
\log \!\frac{\vert I_{i,j,k'+ij} \vert}{\vert I_{i,j,k'}\vert}\right \vert\!.$$
\end{small}The last two terms in this sum are easy to estimate, as the indices $k,k'$ do not mix in none of these. 
Hence, we need to estimate the first term. More precisely, we need to find an upper bound for 
$$\frac{1}{{\vert x-y\vert^\alpha}} \left \vert \log  \frac{|I| |J'|}{|I'||J|}\right \vert,$$
where $I:=I_{i,j,k}$, $J:=f(I)=I_{i,j,k+ij}$ and $I':=I_{i,j,k'}$, $J':=f(I')=I_{i,j,k'+ij}$. 

\vs

To deal with Cases 1 and 2 along the lines explained above, we introduce some notation. 
We say that two real-valued functions $f,g$ satisfy $f\prec g$ if there is a constant $M$ 
such that $|f (x)| \leq M g(x)$ holds for all $x$. Observe that with this notation, for every $a>0$,  
one has $(x+y)^a\prec \max\{|x|^a,|y|^a\}$. When $f$ and $g$ are non-negative functions, we 
will write $f\asymp g$ whenever $f\prec g $ and $g\prec f$. For instance, for $a>0$,  
one has $|x+y|^a \asymp \max\{|x|^a,|y|^a\}$.
 
\vsp

We would like to consider the family of functions $k \mapsto 1 + |i|^p + |j|^q + |k|^r$ together with their second derivatives. 
However,  by (i) we have $r < 2$, hence these functions fail to be twice differentiable. This is why we instead consider the function 
$$\varphi(i,j,\xi):=1+ |i|^p+|j|^q+\theta(\xi),$$ 
where $\theta$ is a fixed  $C^2$ function satisfying $\theta(\xi)=|\xi|^r$ for $|\xi|\geq 1$, and $\theta(0)=0$.
We then define the family of functions $$G_{i,j}(\xi):=\log (\varphi(i,j,\xi)).$$

The following inequality will be of great importance for us: Let $a_1, a_2, a_3$ 
and $b$ be non-negative real numbers such that $a_1/p+a_2/q+a_3/r\leq b$. Then since 
$|i|\leq \varphi^{1/p} (i,j,k)$, $|j|\leq \varphi^{1/q} (i,j,k)$ and $k\leq \varphi^{1/r} (i,j,k)$ 
hold for all integers $i,j,k$, we have 
\begin{equation}\label{eq ineq} |i|^{a_1} \, |j|^{a_2}\, |k|^{a_3}\prec \varphi(i,j,k)^b.\end{equation}

\vsp

We also have the following useful

\begin{lem}\label{lem isla} Let $S := 1+|i|^p+|j|^q$, and  suppose $|\xi-k| \leq S^{1/r}+2|ij|$. 
Then\footnote{Please notice that here (an also below) we are slightly  
abusing of  the notation $\asymp$. Indeed, the precise conclusion should be that there is a universal constant $M$ 
such that $\frac{1}{M} \varphi(i,j,k)\leq \varphi(i,j,\xi)\leq M\varphi(i,j,k)$ holds whenever $|\xi - k| \leq S^{1/r} + 2|ij| $.  } 
$$\varphi(i,j,\xi)\asymp \varphi(i,j,k).$$
\end{lem}

\noindent{\em Proof:} 
By  symmetry, it is enough to find a uniform bound for $\frac{\varphi(i,j,\xi)}{\varphi(i,j,k)}$, and this  follows from
$$\frac{\varphi(i,j,\xi)}{\varphi(i,j,k)} = \frac{\varphi(i,j,k+(\xi-k))}{\varphi(i,j,k)} 
\prec 1 + \frac{|\xi - k|^{r}}{\varphi(i,j,k)}\prec 1+\frac{S+2^r|ij|^{r}}{\varphi(i,j,k)}\leq 2 +\frac{2^r|ij|^{r}}{\varphi(i,j,k)}
$$
and the last expression is  bounded due to conditions (ii) and (iii). 
$\hfill\square$

\vsp\vsp

Now, consider the expression
$$\log \frac{|I| |J'|}{|I'||J|}=\log |I| +\log |J'|-\log|I'|-\log |J|. $$
This can be seen as a ``second increment'' of the function $G_{i,j}$. Indeed, it equals
$$G_{i,j}(k+a+b)-G_{i,j}(k+a)- G_{i,j}(k+b)+ G_{i,j}(k),$$
where, in case 1, $a=-1$ and $b=ij$, and in case 2, $a=k'-k$ and $b=ij$. 
An application of the Mean Value Theorem then yields
\begin{equation}\label{eq referee} G_{i,j}(k+a+b)-G_{i,j}(k+a)- G_{i,j}(k+b)+ G_{i,j}(k)=abG_{i,j}''(\xi)\end{equation}
where $\xi$ is a certain point in $conv\{k,k+a,k+b,k+a+b\}$, the convex hull of $k$, $k+a$,  $k+b$, $k+a+b$.  

Now, for $\xi \notin [-1,1]$, we have that 
\begin{equation}\label{frac-1}
G'_{i,j}(\xi)=\frac{\varphi'}{\varphi}= \pm \frac{r \xi^{r-1}}{\varphi(i,j,\xi)}\prec \frac{ \xi^{r-1}}{\varphi(i,j,\xi)}
\end{equation}
and
\begin{equation}\label{frac-2}
G''_{i,j}(\xi)=\frac{\varphi''}{\varphi}-\left( \frac{\varphi'}{\varphi}\right)^2
=
\pm \frac{r(r-1) \xi^{r-2}}{\varphi(i,j,\xi)} -\frac{r^2 \xi^{2r-2}}{\varphi(i,j,\xi)^2}\prec \frac{\xi^{r-2}}{\varphi(i,j,\xi)},
\end{equation}
where the last bound holds since $\frac{r^2 \xi^{2r-2}}{\varphi(i,j,\xi)^2}= \frac{r^2 \xi^r}{\varphi(i,j,\xi)}\frac{\xi^{r-2}}{\varphi(i,j,\xi)}$, 
and the first factor of this product is always smaller than $r^2$. Besides, for $\xi \in [-1,1]$, the numerators of the right-side 
expressions in (\ref{frac-1}) and (\ref{frac-2}) are bounded from above by some constant independent of $i,j$. 
Therefore, for a general $\xi$, we have $G''_{i,j}(\xi)\prec \frac{1}{\varphi(i,j,\xi)}$.

\vs

Next, consider Case 1, that is, assume that $x$ and $y$ belong to the same interval. By (\ref{eq referee}) and (\ref{frac-2}), we have
$$G_{i,j}(k+ij-1)-G_{i,j}(k-1)- G_{i,j}(k+ij)+ G_{i,j}(k)\prec |i||j|\frac{1}{\varphi(i,j,\xi)}, $$ where $\xi$ is certain point in 
$conv\{ k, k-1, k+ij, k+ij-1\}$. But conditions (iv) and (v) imply that  $|i||j|\prec \varphi^{1-\alpha}(i,j,k)$, and 
since changing $k$ by $k\pm 1$ does not change the asymptotic behavior of $\varphi(i,j,k)$, from Lemma \ref{lem isla}  we have
\begin{equation}\notag \label{eq case 1} \log \frac{|I| |J'|}{|I'||J|}\prec |i||j|\frac{1}{\varphi(i,j,\xi)}\prec \varphi^{-\alpha}(i,j,k).\end{equation}
In particular, except for finitely many indices $(i,j,k)\in \Z^3$, the value of $\varphi(i,j,k)^{\alpha}  \log \frac{|I| |J'|}{|I'||J|}$ is uniformly small. 
Therefore, 
$$\frac{1}{|I|^\alpha} \left| \frac{|I| |J'|}{|I'||J|} -1 \right|\prec  \frac{1}{|I|^\alpha }\log \frac{|I| |J'|}{|I'||J|}\prec 1$$
holds for all indices, as desired.

\vs\vs

Now consider Case 2, namely when $x$ and $y$ belong to different intervals $I := I_{i,j,k}$ and $I' := I_{i,j,k'}$, respectively. 
In this case, by (\ref{eq referee}) and (\ref{frac-2}), we have
\begin{equation}\label{eq-different}
G_{i,j}(k'+ij)- G_{i,j}(k') -G_{i,j}(k+ij)+G_{i,j}(k)\prec \big| ij (k'-k) \big|  \frac{|\xi^{r-2}|}{\varphi(i,j,\xi)},
\end{equation}
where $\xi$ is a certain point in $conv\{k,k',k+ij,k'+ij\}$.
For simplicity, we will assume that $k' - k \geq 2$: the case $k' = k+1$ follows from the previous one using 
property ({\em Derivatives at the endpoints}) just comparing at the right endpoints of $I_{i,j,k}$. Further, we 
also assume $k,k'$ to be positive (the case where both are negative follows by symmetry, and if they have 
different sign, it suffices to consider an intermediate comparison with the term corresponding to $k'' = 0$).  
Finally, by eventually using the triangular inequality, we can restrict ourselves to three different regimes, namely 
when both $k$ and $k'$ belong to $[0, 2|ij|]$, or to $[2|ij|+1,S^{1/r}]$, or to $[S^{1/r}+1,\infty)$, where, as in Lemma 
\ref{lem isla}, we denote $S:=1+|i|^p+|j|^q$. In the same way, we can assume that $x$ is the left end-point of $I$ 
and $y$ is the right end-point of $I'$.

\vsp

Observe that the division into three intervals above is rather natural. The magnitude $S^{1/r}$ marks the point after which the size of the interval 
$I_{i,j,k}$ depends mainly on $k$ and is comparable to $\frac{1}{k^r}$. This will be important to estimate the magnitude $|x-y|$ when $k$ and $k'$ 
are very far apart. 

We start by noticing that in general 
$$\frac{1}{|x-y|^\alpha}=\left(\frac{1}{\sum_{\ell=k}^{k'} |I_{j,k,\ell}|}\right)^\alpha\leq \left( \frac{1}{|k'-k| |I_{i,j,k'}|}\right)^\alpha.$$
Besides, if $k,k'$ are such that $k'-k\leq S^{1/r}$, then $|I|$ and $|I'|$ are comparable by Lemma \ref{lem isla}, hence
\begin{equation} \label{eq isla1}\frac{1}{|x-y|^\alpha}\prec \left( \frac{1}{|k' - k| |I_{i,j,k}|}\right)^\alpha.\end{equation}

We next separately deal with the three regimes $[0, 2|ij|]$, $[2|ij|+1,S^{1/r}]$, and $[S^{1/r}+1,\infty)$.

\begin{itemize}
\item 
Assume that $k$ and $k'$ lie in the interval $[0,2|ij|]$. By Lemma \ref{lem isla}, we have
$$\frac{1}{\varphi(i,j,\xi)}\asymp \frac{1}{\varphi(i,j,k)}.$$
Moreover, from conditions (iv) and (v) we deduce that $|ij||k'-k|\leq 2|ij|^2 \prec \varphi(i,j,k)^{1-\alpha}$. 
Also, as $r < 2$, we have $|\xi^{r-2}|\leq 1$. Therefore, by (\ref{eq-different}), 
$$\frac{1}{{\vert x-y\vert^\alpha}} \left \vert \log  \frac{|I| |J'|}{|I'||J|}\right \vert \prec \frac{1}{|x-y|^\alpha}\frac{|ij||k'-k|}{\varphi(i,j,k)} 
\prec \frac{\varphi^{-\alpha}(i,j,k)}{|x-y|^\alpha} \prec \frac{\varphi^{-\alpha}(i,j,k)}{|I|^\alpha} = 1,$$
as desired.

\vsp

\item 
Assume that $k$ and $k'$ lie in the interval $[2|ij|+1,S^{1/r}]$. Then, 
by (\ref{eq-different}) and  (\ref{eq isla1}), we need to obtain an upper bound for 
$$\left( \frac{1}{|k'-k| |I_{i,j,k}|}\right)^\alpha \big| ij (k'-k) \big|  \frac{|\xi^{r-2}|}{\varphi(i,j,\xi)}.$$
Since $S^{1/r}+|ij|\geq \xi\geq |ij|$ and $r<2$, using Lemma \ref{lem isla} this reduces to estimating the expression 
$$\left( \frac{1}{ |I_{i,j,k}|}\right)^\alpha |ij| \,(k'-k)^{1-\alpha}  \frac{|ij|^{r-2}}{\varphi(i,j,k)}\; =\; |k'-k|^{1-\alpha}  \frac{|ij|^{r-1}}{\varphi(i,j,k)^{1-\alpha}}.$$
In other words, it is enough to show that $(k'-k)^{1-\alpha} |ij|^{r-1} \prec \varphi(i,j,k)^{1-\alpha}$. 
But since $k'-k\leq S^{1/r}$, this reduces to showing that 
\begin{equation}\label{eq isla2}\left( |i|^{p(1-\alpha)/r} + |j|^{q(1-\alpha)/r}\right) |ij|^{r-1} \prec \varphi(i,j,k)^{1-\alpha}. \end{equation}
To show this, we claim that 
$$ |i|^{p(1-\alpha)/r}  |ij|^{r-1} \prec \varphi(i,j,k)^{1-\alpha}$$
follows from (\ref{eq ineq}). (The same inequality changing $|i|^{p(1-\alpha)/r}$ by $|j|^{q(1-\alpha)/r}$ follows in analogous way.) 
Indeed, in order to apply (\ref{eq ineq}) we need to check that 
$$\frac{1-\alpha}{r} + \frac{r-1}{p} + \frac{r-1}{q} \leq 1 - \alpha.$$
However, by (vi), we have $1/p+1/q \leq 1 - \frac{1}{r}$. Therefore, it suffices to show that 
$$\frac{1-\alpha}{r} + (r-1) \big( 1 - \frac{1}{r} \big) \leq 1-\alpha,$$
that is, $\alpha + r \leq 2,$ which is nothing but condition (i).

\vsp

\item Finally, assume that $k$ and $k'$ are in the interval $[S^{1/r},\infty]$. 
If $k'\leq 2k$, then 
$$\frac{\varphi(i,j,k')}{\varphi(i,j,k)} = \frac{\varphi(i,j,k+(k'-k))}{\varphi(i,j,k)} 
\prec 1 + \frac{|k' - k|^{r}}{\varphi(i,j,k)}\leq 1+\frac{|k|^{r}}{\varphi(i,j,k)}\leq 2.$$ 
Therefore, (\ref{eq isla1}) still applies, so that we may proceed as in the second regime case above.  
One then easily checks that, instead of (\ref{eq isla2}), now one needs to show that
$$|k|^{1-\alpha} |ij|^{r-1} \prec \varphi(i,j,k)^{1-\alpha},$$
which still holds thanks to (\ref{eq ineq}) as above.

Assume now that $k'\geq 2k$. The key point in this case is that 
\begin{eqnarray*}
|x-y|
&=&
\sum_{\ell=k}^{k'} |I_{i,j,\ell}|= \sum_{\ell=k}^{k'} \frac{1}{1+|i|^p+|j|^q+|\ell|^r} \\ 
&\succ& \sum_{\ell=k}^{k'} \frac{1}{|\ell|^r}\\ 
&\succ& \int_{k}^{k'} \frac{dx}{x^r}\\
&=& \frac{1}{r-1}\left(\frac{1}{k^{r-1}} -\frac{1}{k'^{r-1}}\right)\\
&\geq& \frac{1}{r-1}\left(1-\frac{1}{2^{r-1}}\right)\frac{1}{k^{r-1}}.
\end{eqnarray*}
Thus, if we further estimate both $\log |I| -\log |J|$ and $\log |I'| -\log |J'|$ 
using the Mean Value Theorem, then we obtain that
$$\frac{1}{{\vert x-y\vert^\alpha}} \left \vert \log  \frac{|I| |J'|}{|I'||J|}\right \vert 
\prec k^{\alpha(r-1)} |ij| \left( \frac{\xi^{r-1}}{\varphi(i,j,\xi)} -\frac{\tilde{\xi}^{r-1}}{\varphi(i,j,\tilde\xi)}\right)$$
for some points $\xi \!\in\! conv\{k,k+ij\}$ and $\tilde \xi \! \in\! conv\{k',k'+ij\}$. 
Since $ \xi \mapsto \frac{\xi }{\varphi(i,j,\xi)}$ is a decreasing function, it suffices to obtain an upper bound for 
$$k^{\alpha(r-1)} |ij|  \frac{k^{r-1}}{\varphi(i,j,k)}.$$
In other words, we need to show that $k^{\alpha(r-1)} |ij| k^{r-1}\prec \varphi(i,j,k)$, which,  by (\ref{eq ineq}), 
follows provided we check that 
$$\frac{(\alpha+1)(r-1)}{r}+1/p+1/q\leq 1.$$
But due to (vi), this holds whenever
$$\frac{(\alpha+1)(r-1)}{r} \leq \frac{1}{r},$$
that is, $(\alpha + 1) (r - 1) \leq 1$, or equivalently, $\alpha r + r \leq \alpha + 2$. However,  
using (i), we obtain
$$\alpha r + r 
= \alpha (r - 1) + (\alpha + r) \leq \alpha + 2,$$
as desired.
\end{itemize}

This finishes the proof of Theorem A.

\vsp\vsp\vsp

\vspace{0.45cm}

\begin{small}

\noindent{\bf Acknowledgments.} We are all grateful to R.Tessera and G.Castro for useful discussions on the subject. Our gratitude also goes to the 
anonymous referee for suggesting us many changes that greatly improved the exposition and correctness of this article. All three authors were funded 
by the Center of Dynamical Systems and Related Fields (Anillo Project 1103 DySyRF,  CONICYT). E.Jorquera was also funded by the Fondecyt 
Project 11121316, and C.Rivas by the Fondecyt Project 1150691 and Inserci\'on 79130017. 


\vsp

\noindent Eduardo Jorquera (eduardo.jorquera@pucv.cl)\\
Instituto de Matem\'atica\\
Pontificia Universidad Cat\'olica de Valpara\'{\i}so\\
Blanco Viel 596, Cerro Bar\'on, Valpara\'{\i}so, Chile

\vsp\vsp\vsp

\noindent Andr\'es Navas (andres.navas@usach.cl) and Crist\'obal Rivas (cristobal.rivas@usach.cl)\\
\noindent Departamento de Matem\'atica y Ciencia de la Computaci\'on\\ 
Universidad de Santiago de Chile\\
Alameda 3363, Estaci\'on Central, Santiago, Chile

\end{small}

\end {document}